\documentclass[fleqn]{article}
\usepackage{amsmath}
\usepackage{amsfonts}
\begin{document}
\renewcommand{\baselinestretch}{1.2}\
\begin{center}
{\Large\bf {Jordan product determined points in matrix algebras }}\vspace{0.4cm}
 \\
{\bf Wenlei Yang \footnote{E-mail address: yang121045760@yahoo.cn}, Jun Zhu{\footnote{E-mail address: zhu$\_$gjun@yahoo.com.cn}}}\\
\vspace{0.3cm} \small{Institute of Mathematics, Hangzhou Dianzi
University, Hangzhou 310018, People's Republic of China}
\end{center}
\underline{~~~~~~~~~~~~~~~~~~~~~~~~~~~~~~~~~~~~~~~~~~~~~~~~~~~~~~~~~~~~~~~~~
~~~~~~~~~~~~~~~~~~~~~~~~~~~~~~~~~~~~~~~~~~~~~~~~~~~~~~~~~~~~~~~}
\vspace{0.1cm}\ \vspace{0.2cm} {\bf{Abstract}}^^L {Let $M_n(R)$ be the algebra of all $n\times n$ matrices over
a unital commutative ring $R$ with 6 invertible. We say that $A\in M_n(R)$ is a Jordan product determined point
if for every $R$-module $X$ and every symmetric $R$-bilinear map $\{\cdot, \cdot\}$ : $M_n(R)\times M_n(R)\to X$ the following two conditions are equivalent: (i) there exists a fixed element $w\in X$ such that $\{x,y\}=w$ whenever $x\circ y=A$, $x,y\in M_n(R)$; (ii) there exists an $R$-linear map $T:M_n(R)^2\to X$ such that $\{x,y\}=T(x\circ y)$ for all $x,y\in M_n(R)$. In this paper, we mainly prove that all the matrix units are the Jordan product determined points in $M_n(R)$ when $n\geq 3$. In addition, we get some corollaries by applying the main results.\vspace{0.2cm}
\\
\vspace{0.2cm} {\it{AMS Classification}}: 15A04
15A27}\\
\vspace{0.2cm}
\underline{~~~~~~~~~~~~~~~~~~~~~~~~~~~~~~~~~~~~~~~~~~~~~~~~~~~~~~~~~~~~~~~~~
~~~~~~~~~~~~~~~~~~~~~~~~~~~~~~~~~~~~~~~~~~~~~~~~~~~~~~~~~~~~~~~}
\vspace{0.1cm}\

\section*{1. Introduction }
~

In this paper, we will mainly discuss Jordan product determined points on matrix algebras. Before proceeding let us fix some symbols and notations in this paper. Let $M_n(R)$ be the algebra of all $n\times n$ matrices over a unital commutative ring $R$ with 6 invertible. Matrix units are denoted by $e_{ij}$ and the Jordan product ``$\circ$'' is defined as $x\circ y=xy+yx$. The identity matrix is denoted by $I$. By $M_n(R)^2$ we denote the $R$-linear span of all elements of the form $xy$ where $x,y\in M_n(R)$.

The concept of zero product (resp. Jordan product, Lie product) determined algebras was introduced by Bre\v{s}ar et al. \cite{bresar}. According to \cite{bresar}, $M_n(B)$ ($n\geq 2$) is zero product determined where $B$ is a unital algebra. If $B$ is a unital algebra with 2 invertible, then $M_n(B)$ ($n\geq 3$) is zero Jordan product determined. From the results above, we can study the linear maps preserving zero product (resp. Jordan product) and derivable (resp. Jordan derivable) at zero point respectively. Wang et al. \cite{Wang1,Wang2} showed that (1) if a symmetric bilinear map $\{\cdot ,\cdot\}$ : $M_n(R)\times M_n(R)\to X$ satisfies the condition that $\{u,u\}=\{e,u\}$ whenever $u^2=u$, then there exists a linear map $f$ from $M_n(R)$ to $X$ such that $\{x,y\}=f(x\circ y)$ for all $x,y\in M_n(R)$; and (2) if an invertible linear map $\delta$ on $M_n(R)$ preserves identity-product, then it is a Jordan automorphism; and a linear map $\sigma$ on $M_n(R)$ is derivable at the identity matrix if and only if it is an inner derivation. Zhu et al. \cite{Zhu1} showed that for every $G\in M_n$, det$G=0$, is an all-multiplicative point in $M_n$. Gong and Zhu \cite{Gong} considered the case of Jordan all-multiplicative point in $M_n$. Zhu et al. \cite{Zhu2} showed that a matrix $G$ is an all-derivable point in $M_n$ if and only if $G\neq 0$. Zhao et al. \cite{Zhao} showed that every element of the algebra of all upper triangular matrices is a Jordan all-derivable point.

Motivated by the concepts and results above, we will consider Jordan product determined points in matrix algebras. For $A\in M_n(R)$, we say that $A$ is a Jordan product determined point if for every $R$-module $X$ and every symmetric $R$-bilinear map $\{\cdot, \cdot\}$ : $M_n(R)\times M_n(R)\to X$ the following two conditions are equivalent: (i) there exists a fixed element $w\in X$ such that $\{x,y\}=w$ whenever $x\circ y=A$, $x,y\in M_n(R)$; (ii) there exists an $R$-linear map $T:M_n(R)^2\to X$ such that $\{x,y\}=T(x\circ y)$ for all $x,y\in M_n(R)$. We say that $G\in M_n(R)$ is a Jordan all-multiplicative point in $M_n(R)$ if for every $M_n(R)$-module $X$ and every Jordan multiplicative $R$-linear map $\varphi$ : $M_n(R)\to X$ at $G$ (i.e. $\varphi (S\circ T)=\varphi (S)\circ \varphi (T)$ for any $S,T\in M_n(R)$, $S\circ T=G$) with $\varphi (I)=I$ is a multiplicative mapping in $M_n(R)$. We say that $H\in M_n(R)$ is a Jordan all-derivable point in $M_n(R)$ if for every $M_n(R)$-module $X$ and every Jordan  derivable $R$-linear map $\varphi$ : $M_n(R)\to X$ at $H$ (i.e. $\varphi (S\circ T)=\varphi (S)\circ T+S\circ \varphi (T)$ for any $S,T\in M_n(R)$, $S\circ T=H$) with $\varphi (I)=0$ is a Jordan derivation in $M_n(R)$. The above two definitions are somewhat different from \cite{Gong} and \cite{Zhao}. In this paper, we will prove that every matrix unit $e_{ij}$ is a Jordan product determined point in $M_n(R)$ when $n\geq 3$. As an application of the result above, we will show that every matrix unit $e_{ij}$ is a Jordan all-multiplicative point and a Jordan all-derivable point respectively.

This paper is organized as follows. Section 2 concerns Jordan product determined points in $M_n(R)$, and we obtain the major results Theorem 2.2 and Theorem 2.3 in this paper. In Section 3, we get some corollaries by applying the main results in section 2.

\section*{2. Jordan product determined points in $M_n(R)$}
~

According to \cite{bresar}, we give the lemma below.\\~\\ {\bf Lemma 2.1.} For $A\in M_n(R)$, $A$ is a Jordan product determined point if and only if for every $R$-module $X$ and every symmetric $R$-bilinear map $\{\cdot ,\cdot \}$ satisfy the condition (i), the following condition (iii) holds true.

 (iii) for every $x_t,y_t\in M_n(R)$ with $\sum_{t=1}^lx_t\circ y_t=0$, $t=1,2,\dots ,l$, $\sum_{t=1}^l\{x_t,y_t\}=0$ holds true.\\~\\ {\bf Proof.} Obviously, the ``only if'' part holds true.\\Conversely, if the condition (iii) holds true, we can define $R$-linear map $T$ : $M_n(R)^2\to X$ as $T(\sum_tx_t\circ y_t)=\sum_t\{x_t,y_t\}$ according to \cite{bresar}. Then $T$ satisfies condition (ii). We only need to prove that $T$ is well defined. Indeed, if condition (iii) is fulfilled, $T$ is well defined obviously. Hence $A$ is a Jordan product determined point. $\Box$\\~\\ {\bf Theorem 2.2.} $e_{ss}$, $s\in \{1,2,\dots ,n\}$, is a Jordan product determined point in $M_n(R)$ when $n\geq 3$.\\~\\ {\bf Proof.} Let $s$ be a fixed number, $X$ be an $R$-module, $\{\cdot,\cdot \}$ : $M_n(R)\times M_n(R)\to X$ be a symmetric (i.e. $\{x,y\}=\{y,x\}$) $R$-bilinear map. Now we assume that there exists a fixed element $w\in X$ such that $\{x,y\}=w$ whenever $x\circ y=e_{ss}$, $x,y\in M_n(R)$. Throughout the proof, $i,j,k,m$ will denote arbitrary indices.

We begin by noticing that $(\frac {1}{2}e_{ss})\circ e_{ss}=e_{ss}$ and as $\{\cdot,\cdot\}$ is symmetric, then
\begin{equation}
\{\frac {1}{2}e_{ss},e_{ss}\}=w=\{e_{ss},\frac {1}{2}e_{ss}\}.
\end{equation}
Next we suppose $n\geq 3$ and divide the proof into three steps.

Step 1. In this step, we assume $i\neq s,j\neq s,k\neq s \; \mathrm{and} \; m\neq s$.

Case 1.1. $i\neq m \; \mathrm{and} \; j\neq k$.\\Since $(\frac {1}{2}e_{ss}+e_{ij})\circ e_{ss}=e_{ss}$ and as $\{\cdot, \cdot\}$ is symmetric, we have that
\begin{equation}
\{e_{ij},e_{ss}\}=0=\{e_{ss},e_{ij}\}.
\end{equation}
Noting $(\frac {1}{2}e_{ss}+e_{ij})\circ (e_{ss}+e_{km})=e_{ss}$, it follows $\{\frac {1}{2}e_{ss}+e_{ij},e_{ss}+e_{km}\}=w$. As $\{\cdot, \cdot\}$ is symmetric, applying (1) and (2), this yields
\begin{equation}
\{e_{ij},e_{km}\}=0.
\end{equation}

Case 1.2. $i,j,k$ are distinct.\\From $(\frac {1}{2}e_{ss}+e_{ik})\circ (e_{ss}+e_{kk}-e_{ii})=e_{ss}$, we obtain that $\{\frac {1}{2}e_{ss}+e_{ik}, e_{ss}+e_{kk}-e_{ii}\}=w$. Because $\{\cdot, \cdot\}$ is symmetric, it follows from (1) and (2) that
\begin{equation}
\{e_{ik},e_{kk}\}=\{e_{ik},e_{ii}\}=\{e_{ii},e_{ik}\}=\{e_{kk},e_{ik}\}.
\end{equation}
Now we assume $n>3$. As $(\frac {1}{2}e_{ss}+e_{ij}+e_{ik})\circ (e_{ss}+e_{jk}-e_{kk})=e_{ss}$, we derive $\{\frac {1}{2}e_{ss}+e_{ij}+e_{ik}, e_{ss}+e_{jk}-e_{kk}\}=w$. Using (1), (2) and (3), this can be reduced to
\begin{equation}
\{e_{ij},e_{jk}\}=\{e_{ik},e_{kk}\}.
\end{equation}
(4) together with (5) yield
\begin{equation}
\{e_{ij}, e_{jk}\}=\{e_{ik}, e_{kk}\}=\{e_{ik}, e_{ii}\}=\{e_{ii}, e_{ik}\}.
\end{equation}
Since $\{\cdot, \cdot\}$ is symmetric, it follows that
\begin{equation}
\{e_{ii}, e_{ik}\}=\{e_{ik}, e_{ii}\}=\{e_{kk}, e_{ik}\}=\{e_{jk},e_{ij}\}.
\end{equation}

Case 1.3. $i\neq j$.\\By $(\frac{1}{2}e_{ss}+\frac{1}{2}e_{ii}+e_{ij}-\frac{1}{2}e_{jj})\circ (e_{ss}+e_{ji}-e_{ii}+e_{jj})=e_{ss}$, it is clear that $\{\frac {1}{2}e_{ss}+\frac {1}{2}e_{ii}+e_{ij}-\frac {1}{2}e_{jj}, e_{ss}+e_{ji}-e_{ii}+e_{jj}\}=w$. Then we can get from (1), (2), (3) and (4) that
\begin{equation}
\{e_{ij}, e_{ji}\}=\frac{1}{2}\{e_{ii}, e_{ii}\}+\frac{1}{2}\{e_{jj}, e_{jj}\}.
\end{equation}

Step 2. In this step, we consider some of the indices of the matrix units equal to $s$.

Case 2.1. $i\neq s, j\neq s, k\neq s \; \mathrm{and} \; m\neq s$.\\For $(\frac{1}{2}e_{ss}-e_{si})\circ (e_{ss}-e_{ii})=e_{ss}$, we have that $\{\frac {1}{2}e_{ss}-e_{si}, e_{ss}-e_{ii}\}=w$. Applying (1) and (2), it can be reduced to $\{e_{si}, e_{ss}\}=\{e_{si}, e_{ii}\}$. As $\{\cdot, \cdot\}$ is symmetric, then it follows that
\begin{equation}
\{e_{si}, e_{ss}\}=\{e_{si}, e_{ii}\}=\{e_{ss}, e_{si}\}=\{e_{ii}, e_{si}\}.
\end{equation}
Since $(\frac{1}{2}e_{ss}-e_{is})\circ (e_{ss}-e_{ii})=e_{ss}$, a similar discussion as above shows that \begin{equation}
\{e_{is}, e_{ss}\}=\{e_{is}, e_{ii}\}=\{e_{ss}, e_{is}\}=\{e_{ii}, e_{is}\}.
\end{equation}
Let $j\neq k$, then we have $(\frac {1}{2}e_{ss}+e_{sj})\circ (e_{ss}+e_{km}-e_{jj})=e_{ss}$ and it implies that $\{\frac {1}{2}e_{ss}+e_{sj}, e_{ss}+e_{km}-e_{jj}\}=w$. As $\{\cdot, \cdot\}$ is symmetric and by (1), (2) and (9), this yields
\begin{equation}
\{e_{sj}, e_{km}\}=0=\{e_{km}, e_{sj}\}\;\mathrm {if}\; j\neq k.
\end{equation}
From $(\frac{1}{2}e_{ss}+e_{si})\circ (e_{ss}-2e_{si})=e_{ss}$, we derive $\{\frac{1}{2}e_{ss}+e_{si}, e_{ss}-2e_{si}\}=w$. Then it follows from (1) and (9) that
\begin{equation}
\{e_{si}, e_{si}\}=0.
\end{equation}
Since $(\frac {1}{2}e_{ss}+e_{sj}-\frac{1}{2}e_{si})\circ (e_{ss}+e_{si}-e_{jj})=e_{ss}$ when $i\neq j$, then we have $\{\frac {1}{2}e_{ss}+e_{sj}-\frac{1}{2}e_{si}, e_{ss}+e_{si}-e_{jj}\}=w$. Using (1), (2), (9), (11) and (12), it is clear that
\begin{equation}
\{e_{sj}, e_{si}\}=0\;\mathrm {if}\; i\neq j.
\end{equation}
If $j\neq m$, from $(\frac {1}{2}e_{ss}+e_{js})\circ (e_{ss}+e_{km}-e_{jj})=e_{ss}$ we can get that $\{\frac {1}{2}e_{ss}+e_{js}, e_{ss}+e_{km}-e_{jj}\}=w$. As $\{\cdot, \cdot\}$ is symmetric and by (1), (2) and (10), it follows
\begin{equation}
\{e_{js}, e_{km}\}=0=\{e_{km}, e_{js}\}\;\mathrm {if}\; j\neq m.
\end{equation}
For $(\frac{1}{2}e_{ss}+e_{is})\circ (e_{ss}-2e_{is})=e_{ss}$, we have that $\{\frac{1}{2}e_{ss}+e_{is}, e_{ss}-2e_{is}\}=w$. Applying (1) and (10), it can be reduced to
\begin{equation}
\{e_{is}, e_{is}\}=0.
\end{equation}
Assume $i\neq j$, it follows from $(\frac {1}{2}e_{ss}+e_{js}-\frac{1}{2}e_{is})\circ (e_{ss}+e_{is}-e_{jj})=e_{ss}$  that $\{\frac {1}{2}e_{ss}+e_{js}-\frac{1}{2}e_{is}, e_{ss}+e_{is}-e_{jj}\}=w$. Using (1), (2), (10), (14) and (15), this yields
\begin{equation}
\{e_{js}, e_{is}\}=0 \;\mathrm {if}\; i\neq j.
\end{equation}

Case 2.2. $i\neq s, k\neq s \; \mathrm{and} \;j\neq s$.\\Since $(\frac {1}{2}e_{ss}+e_{ki}-e_{ks}-\frac {1}{2}e_{ii})\circ (e_{ss}+e_{is})=e_{ss}$ if $k\neq i$, then we have that $\{\frac {1}{2}e_{ss}+e_{ki}-e_{ks}-\frac {1}{2}e_{ii}, e_{ss}+e_{is}\}=w$. By (1), (2), (10) and (16), this yields $\{e_{ki}, e_{is}\}=\{e_{ks}, e_{ss}\}$ if $k\neq i$. For our purpose, it is more convenient to rewrite this equation as $\{e_{ik}, e_{ks}\}=\{e_{is}, e_{ss}\}$ if $i\neq k$. As $\{\cdot, \cdot\}$ is symmetric, then we can conclude from the above equation and (10) that
\begin{equation}
\{e_{ik}, e_{ks}\}=\{e_{is}, e_{ss}\}=\{e_{ii}, e_{is}\}=\{e_{is}, e_{ii}\}=\{e_{ss}, e_{is}\}=\{e_{ks}, e_{ik}\}\;\mathrm {if}\; i\neq k.
\end{equation}
If $i\neq k$, then we have $(\frac {1}{2}e_{ss}+e_{ik}-e_{sk}-\frac{1}{2}e_{ii})\circ (e_{ss}+e_{si})=e_{ss}$. By a similar discussion as above, this yields
\begin{equation}
\{e_{sk}, e_{ki}\}=\{e_{ss}, e_{si}\}=\{e_{si}, e_{ii}\}=\{e_{ii}, e_{si}\}=\{e_{si}, e_{ss}\}=\{e_{ki}, e_{sk}\}\;\mathrm {if}\; i\neq k.
\end{equation}
If $j\neq k$, from $(\frac {1}{2}e_{ss}+e_{sk}+e_{jk}-\frac {1}{2}e_{js})\circ (e_{ss}+e_{js}-e_{kk})=e_{ss}$ we have that $\{\frac {1}{2}e_{ss}+e_{sk}+e_{jk}-\frac {1}{2}e_{js}, e_{ss}+e_{js}-e_{kk}\}=w$. For $\{\cdot, \cdot\}$ is symmetric, applying (1), (2), (14), (15), (17) and (18), it follows
\begin{equation}
\{e_{sk}, e_{js}\}=\{e_{jk}, e_{kk}\}=\{e_{js}, e_{sk}\}\;\mathrm {if}\;j\neq k.
\end{equation}

Case 2.3. $i\neq s$.\\Since $(2e_{si}+\frac {1}{2}e_{is}-e_{ii})\circ (e_{si}+\frac {1}{4}e_{is}+\frac {1}{2}e_{ii})=e_{ss}$, we obtain $\{2e_{si}+\frac {1}{2}e_{is}-e_{ii}, e_{si}+\frac {1}{4}e_{is}+\frac {1}{2}e_{ii}\}=w$. Using (1), (12), (15), (17) and (18), this can be reduced to $\{e_{si}, e_{is}\}+\{e_{is}, e_{si}\}=\{e_{ss}, e_{ss}\}+\{e_{ii}, e_{ii}\}$. As $\{e_{si}, e_{is}\}=\{e_{is}, e_{si}\}$, it leads to
\begin{equation}
\{e_{si}, e_{is}\}=\{e_{is}, e_{si}\}=\frac {1}{2}\{e_{ss}, e_{ss}\}+\frac {1}{2}\{e_{ii}, e_{ii}\}.
\end{equation}

Step 3. Now concluding from case 1.1 and case 2.1, we can obtain
\begin{equation}
\{e_{ij}, e_{km}\}=0\;\mathrm{for \;every}\; i, j, k, m,\;\mathrm{if}\; i\neq m \; \mathrm{and} \; j\neq k.
\end{equation}
If $i,j,k$ are distinct and $n=3$, it follows from (4), (17), (18) and (19) that
\begin{equation}
\{e_{ij}, e_{jk}\}=\{e_{ik}, e_{kk}\}=\{e_{ii}, e_{ik}\}=\{e_{ik}, e_{ii}\}=\{e_{kk}, e_{ik}\}=\{e_{jk}, e_{ij}\}.
\end{equation}
If $i,j,k$ are distinct and $n>3$, from (6), (7), (17), (18) and (19) we have the equations above as well. So (22) holds true whenever $n \geq 3$.\\
From case 1.3 and case 2.3, we have
\begin{equation}
\{e_{ij}, e_{ji}\}=\frac {1}{2}\{e_{ii}, e_{ii}\}+\frac {1}{2}\{e_{jj}, e_{jj}\}\;\mathrm {if}\;i\neq j.
\end{equation}
Let $\sum_{t=1}^lx_t\circ y_t=0$ where $x_t, y_t\in M_n(R)$, $t=1,2,\ldots ,l$. We write $x_t$ and $y_t$ for
\begin{equation}\nonumber
x_t=\sum_{i=1}^n\sum_{j=1}^na_{ij}^te_{ij},\; y_t=\sum_{k=1}^n\sum_{m=1}^nb_{km}^te_{km}.
\end{equation}
Then for all $i$ and $m$ we have that
\begin{equation}
\sum_{t=1}^l\sum_{j=1}^n(a_{ij}^tb_{jm}^t+b_{ij}^ta_{jm}^t)=0.
\end{equation}
Now we will show $\sum_{t=1}^l\{x_t, y_t\}=0$. Note that
\begin{eqnarray*}
\sum_{t=1}^l\{x_t, y_t\}&=&\sum_{t=1}^l\sum_{i=1}^n\sum_{j=1}^n\sum_{k=1}^n\sum_{m=1}^n\{a_{ij}^te_{ij}, b_{km}^te_{km}\}\\&=&\sum_{t=1}^l\sum_{i=1}^n\sum_{j=1}^n\sum_{k=1}^n\sum_{m=1}^na_{ij}^tb_{km}^t\{e_{ij}, e_{km}\}.
\end{eqnarray*}
According to our assumptions, it follows from (21), (22), (23) and (24) that
\begin{eqnarray*}
\sum_{t=1}^l\{x_t, y_t\} &= &\sum_{t=1}^l\sum_{i=1}^n\sum_{j=1}^n\sum_{k=1 \atop k\neq j}^n\sum_{m=1 \atop m\neq i}^na_{ij}^tb_{km}^t\{e_{ij}, e_{km}\}+\sum_{t=1}^l\sum_{i=1}^n\sum_{j=1}^n\sum_{m=1 \atop m\neq i}^na_{ij}^tb_{jm}^t\{e_{ij}, e_{jm}\}\\
                          &+& \sum_{t=1}^l\sum_{i=1}^n\sum_{j=1}^n\sum_{k=1 \atop k\neq j}^na_{ij}^tb_{ki}^t\{e_{ij}, e_{ki}\}+\sum_{t=1}^l\sum_{i=1}^n\sum_{j=1}^na_{ij}^tb_{ji}^t\{e_{ij}, e_{ji}\}\\
                          &= &\sum_{t=1}^l\sum_{i=1}^n\sum_{j=1}^n\sum_{m=1 \atop m\neq i}^n(a_{ij}^tb_{jm}^t+b_{ij}^ta_{jm}^t)\{e_{ij}, e_{jm}\}\\
                          &+&\frac {1}{2}\sum_{t=1}^l\sum_{i=1}^n\sum_{j=1}^na_{ij}^tb_{ji}^t\left(\{e_{ii}, e_{ii}\}+\{e_{jj}, e_{jj}\}\right)\\
                          &=& \sum_{i=1}^n\sum_{m=1 \atop m\neq i}^n\left[\sum_{t=1}^l\sum_{j=1}^n(a_{ij}^tb_{jm}^t+b_{ij}^ta_{jm}^t)\{e_{im}, e_{mm}\}\right]\\
                          &+&\frac{1}{2}\sum_{i=1}^n\left[\sum_{t=1}^l\sum_{j=1}^n(a_{ij}^tb_{ji}^t+ b_{ij}^ta_{ji}^t)\{e_{ii}, e_{ii}\}\right]\\
                          &=& 0.
\end{eqnarray*}
By Lemma 2.1, $e_{ss}$ is a Jordan product determined point. $\Box$\\~\\ {\bf Theorem 2.3.} $e_{pq}$, $p\neq q$, is a Jordan product determined point in $M_n(R)$ when $n\geq 3$.\\~\\ {\bf Proof.} Let $p, q$ be the distinct fixed indices, $\{\cdot, \cdot\}$ : $M_n(R)\times M_n(R)\to X$ be a symmetric (i.e. $\{x, y\}=\{y, x\}$) $R$-bilinear map where $X$ is an $R$-module. And we assume that there exists a fixed element $w\in X$ such that $\{x,y\}=w$ whenever $x\circ y=e_{pq}$, $x,y\in M_n(R)$. Throughout the proof, $i,j,k,m$ will denote arbitrary indices.

According to Lemma 2.1 and Step 3 in the proof of Theorem 2.2, we only need to verify (21), (22) and (23) hold true when $n\geq 3$. Now we suppose $n\geq 3$ and divide the proof into several steps.

Step 1. In this step, we assume $i\neq q \; \mathrm {and} \; m\neq p$.\\Since $e_{ps}\circ e_{sq}=e_{pq}\:(s=1,2,\ldots ,n)$ and as $\{\cdot, \cdot\}$ is symmetric, hence we have
\begin{equation}
\{e_{ps}, e_{sq}\}=w=\{e_{sq}, e_{ps}\},\; s=1,2,\ldots ,n.
\end{equation}
Choosing $s\neq j\; \mathrm {and} \;s\neq k$, then we obtain $e_{ps}\circ (e_{sq}+e_{km})=e_{pq}$ and $(e_{ps}+e_{ij})\circ e_{sq}=e_{pq}$ respectively. For $\{\cdot, \cdot\}$ is symmetric, applying (25), it follows
\begin{equation}
\{e_{ps}, e_{km}\}=\{e_{km}, e_{ps}\}=\{e_{ij}, e_{sq}\}=\{e_{sq}, e_{ij}\}=0 \;\mathrm{if}\:s\neq j\; \mathrm {and} \; s\neq k.
\end{equation}
Given $s\neq j \; \mathrm {and} \; s\neq k$, then if $i\neq m\; \mathrm {and} \;j\neq k$ we have $(e_{ps}+e_{ij})\circ (e_{sq}+e_{km})=e_{pq}$. It is clear that $\{e_{ps}+e_{ij},e_{sq}+e_{km}\}=w$. As $\{\cdot, \cdot\}$ is symmetric, using (25) and (26), this yields
\begin{equation}
\{e_{ij}, e_{km}\}=0=\{e_{km}, e_{ij}\}\:\mathrm{if}\:i\neq m\; \mathrm {and} \;j\neq k.
\end{equation}

Step 2. In this step, we suppose $i\neq q\; \mathrm {and} \;m\neq p$.\\Assuming $s\neq i\; \mathrm {and} \;s\neq m$, then if $i\neq m\; \mathrm {and} \;i\neq p$ we can verify $(e_{ps}+e_{im})\circ (e_{sq}+e_{mm}-e_{ii})=e_{pq}$. It follows that $\{e_{ps}+e_{im}, e_{sq}+e_{mm}-e_{ii}\}=w$. Since $\{\cdot, \cdot\}$ is symmetric, by (25) and (27), this yields
\begin{equation}
\{e_{im}, e_{mm}\}=\{e_{im}, e_{ii}\}=\{e_{ii}, e_{im}\}\:\mathrm {if}\: i\neq m\; \mathrm {and} \;i\neq p.
\end{equation}
From $e_{pm}\circ (e_{mq}+e_{mm}-e_{pp})=e_{pq}$, we have that $\{e_{pm}, e_{mq}+e_{mm}-e_{pp}\}=w$. As $\{\cdot, \cdot\}$ is symmetric, it follows from (25) that
\begin{equation}
\{e_{pm}, e_{mm}\}=\{e_{pm}, e_{pp}\}=\{e_{pp}, e_{pm}\}.
\end{equation}
Choosing $s\neq j \; \mathrm {and} \; s\neq m$, if $i,j,k$ are distinct we obtain $(e_{ps}+e_{ij}+e_{im})\circ (e_{sq}+e_{jm}-e_{mm})=e_{pq}$. Then it leads to $\{e_{ps}+e_{ij}+e_{im}, e_{sq}+e_{jm}-e_{mm}\}=w$. Applying (25) and (27), this yields
\begin{equation}
\{e_{ij}, e_{jm}\}=\{e_{im}, e_{mm}\}\;\mathrm {if}\;i,j,k\;\mathrm{are distinct}.
\end{equation}
Since $\{\cdot, \cdot\}$ is symmetric, if $i,j,m$ are distinct, we can conclude from (28), (29) and (30) that
\begin{equation}
\{e_{ij}, e_{jm}\}=\{e_{im}, e_{mm}\}=\{e_{ii}, e_{im}\}=\{e_{im}, e_{ii}\}=\{e_{mm}, e_{im}\}=\{e_{jm}, e_{ij}\}.
\end{equation}

Step 3. In this step, we assume $i\neq q\;\mathrm {and}\;m\neq p$.\\Choosing $s\neq i\; \mathrm {and} \; s\neq m$, we suppose $n>3, \:i\neq p,\: i\neq m\; \mathrm {and} \;m\neq q$. As $(e_{ps}+\frac {1}{2}e_{ii}+e_{im}-\frac{1}{2}e_{mm})\circ (e_{sq}+e_{mi}-e_{ii}+e_{mm})=e_{pq}$, it follows that $\{e_{ps}+\frac {1}{2}e_{ii}+e_{im}-\frac{1}{2}e_{mm}, e_{sq}+e_{mi}-e_{ii}+e_{mm}\}=w$. Using (25), (27) and (31), if $n>3,\:i\neq p,\: i\neq m\; \mathrm {and} \;m\neq q$ we have that
\begin{equation}
\{e_{im},e_{mi}\}=\frac{1}{2}\{e_{ii}, e_{ii}\}+\frac{1}{2}\{e_{mm}, e_{mm}\}.
\end{equation}

Step 4. In this step, we assume $i=q\; \mathrm {and} \;m\neq p$.

Case 4.1. $j\neq p\; \mathrm {and} \;k\neq q$.\\By (27), we have $\{e_{km}, e_{qj}\}=0$ if $m\neq q\; \mathrm {and} \;j\neq k$. As $\{\cdot, \cdot\}$ is symmetric, this yields
\begin{equation}
\{e_{qj},e_{km}\}=\{e_{km}, e_{qj}\}=0\;\mathrm {if}\; m\neq q\; \mathrm {and} \;j\neq k.
\end{equation}

Case 4.2. $j\neq p\; \mathrm {and} \; k=q$.\\Since $(e_{pq}+e_{qj})\circ (e_{qq}-e_{jj})=e_{pq}$ if $j\neq q$, then we have $\{e_{pq}+e_{qj}, e_{qq}-e_{jj}\}=w$. For $\{\cdot, \cdot\}$ is symmetric, it follows from (25) and (27) that
\begin{equation}
\{e_{qj}, e_{qq}\}=\{e_{qj}, e_{jj}\}=\{e_{jj}, e_{qj}\}=\{e_{qq}, e_{qj}\}\;\mathrm {if}\; j\neq q.
\end{equation}
Noting $(e_{pq}+e_{qj}+e_{pp})\circ (e_{qq}+e_{qm}-e_{jj}-e_{pm})=e_{pq}$ if $j\neq q\; \mathrm {and} \;m\neq q$, then we obtain $\{e_{pq}+e_{qj}+e_{pp}, e_{qq}+e_{qm}-e_{jj}-e_{pm}\}=w$. Applying (25), (27), (31), (33) and (34), it can be reduced to
\begin{equation}
\{e_{qj}, e_{qm}\}=0 \;\mathrm{if}\; j\neq q\; \mathrm {and} \;m\neq q.
\end{equation}

Case 4.3. $j=p\; \mathrm {and} \;k\neq q$.\\From $(\frac {1}{2}e_{pp}+\frac{1}{2}e_{pq}+e_{qp})\circ (e_{pq}-2e_{qp}+e_{qq})=e_{pq}$, we have that $\{\frac {1}{2}e_{pp}+\frac{1}{2}e_{pq}+e_{qp}, e_{pq}-2e_{qp}+e_{qq}\}=w$. As $\{\cdot, \cdot\}$ is symmetric, it follows from (25) and (27) that
\begin{equation}
\{e_{qp}, e_{qq}\}-\{e_{pp}, e_{qp}\}-2\{e_{qp}, e_{qp}\}=0.
\end{equation}
Since $(\frac{1}{2}e_{pp}+e_{pq}+2e_{qp})\circ (e_{pq}-2e_{qp}+\frac{1}{2}e_{qq})=e_{pq}$, a similar discussion shows that
\begin{equation}
\{e_{qp}, e_{qq}\}-\{e_{pp}, e_{qp}\}-4\{e_{qp}, e_{qp}\}=0.
\end{equation}
Because $\{\cdot, \cdot\}$ is symmetric, comparing (36) with (37), we get that
\begin{equation}
\{e_{qp}, e_{qp}\}=0,
\end{equation}
\begin{equation}
\{e_{qp}, e_{qq}\}=\{e_{pp}, e_{qp}\}=\{e_{qp},e_{pp}\}=\{e_{qq}, e_{qp}\}.
\end{equation}
Noting $(\frac{1}{2} e_{pp}+\frac{1}{2}e_{pq}+e_{qp})\circ (e_{pq}-2e_{qp}+e_{qq}+e_{km})=e_{pq}$ if $k\neq p\; \mathrm {and} \;m\neq q$, we have that $\{\frac{1}{2} e_{pp}+\frac{1}{2}e_{pq}+e_{qp}, e_{pq}-2e_{qp}+e_{qq}+e_{km}\}=w$. As $\{\cdot, \cdot\}$ is symmetric, applying (25), (27) and (36), if $k\neq p\; \mathrm {and} \; m\neq q$ it leads to
\begin{equation}
\{e_{qp}, e_{km}\}=0=\{e_{km}, e_{qp}\}.
\end{equation}

Case 4.4. $j=p\; \mathrm {and} \;k=q$.\\If $m\neq q$, then we can verify that
\begin{equation}\nonumber
(\frac{1}{2}e_{pp}+\frac{1}{2}e_{pq}+e_{qp}-e_{mm})\circ (e_{pq}-2e_{qp}+e_{qq}+e_{qm}+e_{pm})=e_{pq},
\end{equation}
\begin{equation}\nonumber
(\frac{1}{2}e_{pp}+\frac{1}{2}e_{pq}+e_{qp}+\frac {1}{2}e_{mm})\circ (e_{pq}-2e_{qp}+e_{qq}-2e_{qm}+e_{pm})=e_{pq}.
\end{equation}
Hence we have that
\begin{equation}\nonumber
\{\frac{1}{2}e_{pp}+\frac{1}{2}e_{pq}+e_{qp}-e_{mm}, e_{pq}-2e_{qp}+e_{qq}+e_{qm}+e_{pm}\}=w,
\end{equation}
\begin{equation}\nonumber
\{\frac{1}{2}e_{pp}+\frac{1}{2}e_{pq}+e_{qp}+\frac {1}{2}e_{mm}, e_{pq}-2e_{qp}+e_{qq}-2e_{qm}+e_{pm}\}=w.
\end{equation}
For $\{\cdot, \cdot\}$ is symmetric, it follows from (25), (27), (31), (38), (39) and (40) that
\begin{equation}\nonumber
\{e_{qp}, e_{pm}\}+\{e_{qp},e_{qm}\}-\{e_{mm}, e_{qm}\}=0,
\end{equation}
\begin{equation}\nonumber
\{e_{qp}, e_{pm}\}-2\{e_{qp},e_{qm}\}-\{e_{mm}, e_{qm}\}=0.
\end{equation}
Thus if $m\neq q$ we have that
\begin{equation}
\{e_{qp},e_{qm}\}=0=\{e_{qm},e_{qp}\},
\end{equation}
\begin{equation}
\{e_{qp}, e_{pm}\}=\{e_{mm}, e_{qm}\}=\{e_{qm}, e_{mm}\}=\{e_{pm}, e_{qp}\}.
\end{equation}

Step 5. In this step, we assume $i\neq q\; \mathrm {and} \;m=p$.

Case 5.1. $j\neq p\; \mathrm {and} \;k\neq q$.\\From (27), if $i\neq p\; \mathrm {and} \;j\neq k$ it follows that
\begin{equation}
\{e_{ij}, e_{kp}\}=0=\{e_{kp}, e_{ij}\}.
\end{equation}

Case 5.2. $j=p\; \mathrm {and} \;k\neq q$.\\Since $(e_{pq}+e_{ip}-e_{kq})\circ (e_{qq}+e_{kp})=e_{pq}$ if $i\neq p\; \mathrm {and} \;k\neq p$, it is clear that $\{e_{pq}+e_{ip}-e_{kq}, e_{qq}+e_{kp}\}=w$. Applying (25), (27) and (31), it can be reduced to
\begin{equation}
\{e_{ip}, e_{kp}\}=0.
\end{equation}

Case 5.3. $j\neq p\; \mathrm {and} \;k=q$.\\By (40), if $i\neq p\; \mathrm {and} \;j\neq q$ we have that
\begin{equation}
\{e_{ij}, e_{qp}\}=0.
\end{equation}

Case 5.4. $j=p\; \mathrm {and} \;k=q$.\\If $i\neq p\; \mathrm {and} \;i\neq q$, then we can verify that
\begin{eqnarray*}
(\frac{1}{2}e_{pp}+\frac{1}{2}e_{pq}+e_{qp}-e_{ii})\circ (e_{pq}-2e_{qp}+e_{qq}+2e_{ip}+e_{iq})=e_{pq},\\
(\frac{1}{2}e_{pp}+\frac{1}{2}e_{pq}+e_{qp}+\frac{1}{2}e_{ii})\circ (e_{pq}-2e_{qp}+e_{qq}-e_{ip}+e_{iq})=e_{pq}.
\end{eqnarray*}
According to our assumptions, using (25), (27), (31), (38), (39) and (40), we have
\begin{equation}
\{e_{pp},e_{ip}\}+\{e_{qp}, e_{iq}\}-2\{e_{ii}, e_{ip}\}+2\{e_{qp},e_{ip}\}=0,
\end{equation}
\begin{equation}
\{e_{pp},e_{ip}\}-2\{e_{qp}, e_{iq}\}+\{e_{ii}, e_{ip}\}+2\{e_{qp},e_{ip}\}=0.
\end{equation}
Comparing the two equations above, as $\{\cdot, \cdot\}$ is symmetric, it follows that
\begin{equation}
\{e_{qp}, e_{iq}\}=\{e_{ii},e_{ip}\}=\{e_{iq},e_{qp}\}=\{e_{ip},e_{ii}\} \;\mathrm{if}\; i\neq p\; \mathrm {and} \;i\neq q.
\end{equation}
Substituting (48) into (46), we obtain
\begin{equation}
\{e_{pp},e_{ip}\}-\{e_{ii}, e_{ip}\}+2\{e_{qp},e_{ip}\}=0\;\mathrm {if} \;i\neq p\; \mathrm {and} \;i\neq q.
\end{equation}
Noting $(e_{pq}+e_{ip})\circ (e_{ii}-e_{pp}+2e_{qq})=e_{pq}$ if $i\neq p\; \mathrm{and}\;i\neq q$, then we have that $\{e_{pq}+e_{ip}, e_{ii}-e_{pp}+2e_{qq}\}=w$. Since $\{\cdot, \cdot\}$ is symmetric, it follows from (25) and (27) that
\begin{equation}
\{e_{ip}, e_{ii}\}=\{e_{ip}, e_{pp}\}=\{e_{pp}, e_{ip}\}=\{e_{ii}, e_{ip}\}\;\mathrm{if}\; i\neq p \; \mathrm {and} \;i\neq q.
\end{equation}
Substituting (50) into (49), we derive
\begin{equation}
\{e_{qp}, e_{ip}\}=0=\{e_{ip}, e_{qp}\}\;\mathrm{if}\;i\neq p\; \mathrm {and} \;i\neq q.
\end{equation}

Step 6. In this step, we assume $i=q\; \mathrm {and} \;m=p$.

Case 6.1. $j\neq p\; \mathrm {and} \;k\neq q$.\\From (27), if $j\neq k$ it is clear that
\begin{equation}
\{e_{qj},e_{kp}\}=\{e_{kp}, e_{qj}\}=0.
\end{equation}

Case 6.2. $j\neq p\; \mathrm {and} \; k=q$.\\From (41), we have
\begin{equation}
\{e_{qj}, e_{qp}\}=\{e_{qp}, e_{qj}\}=0\;\mathrm{if}\;j\neq q.
\end{equation}

Case 6.3. $j=p\; \mathrm {and} \;k\neq q$.\\From (51), it follows
\begin{equation}
\{e_{qp}, e_{kp}\}=\{e_{kp}, e_{qp}\}=0 \;\mathrm{if}\;k\neq p.
\end{equation}

Case 6.4. $j=p\; \mathrm {and} \;k=q$.\\From (38), we obtain that
\begin{equation}
\{e_{qp}, e_{qp}\}=0.
\end{equation}

Step 7. In this step, we assume $i=q\; \mathrm {and} \;m\neq p$.\\If $n>3,\:j\neq p,\: j\neq q, \:j\neq m\; \mathrm {and} \;m\neq q$, choosing $s\neq q\; \mathrm {and} \;s\neq m$, then we can verify $(e_{ps}+e_{mm}-e_{jm})\circ (e_{sq}+e_{qm}+e_{qj})=e_{pq}$. It follows that $\{e_{ps}+e_{mm}-e_{jm}, e_{sq}+e_{qm}+e_{qj}\}=w$. As $\{\cdot, \cdot\}$ is symmetric, by (25) and (27), if $n>3,\:j\neq p, \:j\neq q,\: j\neq m\; \mathrm {and} \;m\neq q$ we have that
\begin{equation}
\{e_{mm}, e_{qm}\}=\{e_{jm}, e_{qj}\}=\{e_{qj}, e_{jm}\}=\{e_{qm}, e_{mm}\}.
\end{equation}
In the case of $n=3$, if $q, j, m$ are distinct we can conclude from (34) and (42) that
\begin{equation}
\{e_{qj}, e_{jm}\}=\{e_{qm}, e_{mm}\}=\{e_{qq}, e_{qm}\}=\{e_{qm}, e_{qq}\}=\{e_{mm}, e_{qm}\}=\{e_{jm}, e_{qj}\}.
\end{equation}
In the case of $n>3$, if $q, j, m$ are distinct, from (34), (42) and (56) we have the same equations as well. So (57) holds true whenever $n\geq 3$.

Step 8. In this step, we assume that $i\neq q\; \mathrm {and} \; m=p$.\\If $n>3,\:i\neq p,\: i\neq j,\: j\neq p\; \mathrm {and} \; j\neq q$, choosing $s\neq p\; \mathrm {and} \;s\neq i$, then we have that $(e_{ps}+e_{jp}-e_{ip})\circ (e_{sq}+e_{ij}+e_{ii})=e_{pq}$. From the assumptions this yields $\{e_{ps}+e_{jp}-e_{ip}, e_{sq}+e_{ij}+e_{ii}\}=w$. Because $\{\cdot, \cdot\}$ is symmetric, if $n>3,\:i\neq p,\: i\neq j,\: j\neq p\; \mathrm {and} \;j\neq q$ it follows from (25) and (27) that
\begin{equation}
\{e_{jp}, e_{ij}\}=\{e_{ip}, e_{ii}\}=\{e_{ij}, e_{jp}\}=\{e_{ii}, e_{ip}\}.
\end{equation}
If $i,j,p$ are distinct and $n=3$, then we can conclude from (48) and (50) that
\begin{equation}
\{e_{ij}, e_{jp}\}=\{e_{ip}, e_{pp}\}=\{e_{ii}, e_{ip}\}=\{e_{ip}, e_{ii}\}=\{e_{pp}, e_{ip}\}=\{e_{jp}, e_{ij}\}.
\end{equation}
If $i,j,p$ are distinct and $n>3$, from (48), (50) and (58) we obtain the equations above as well. Hence, (59) holds true whenever $n\geq 3$.

Step 9. Selecting $s\neq p\; \mathrm {and} \;s\neq m$, then if $m\neq p\;\mathrm{and} \;m\neq q$ we can verify $(e_{ps}+\frac{1}{2}e_{pp}+e_{mp}-\frac{1}{2}e_{mm})\circ (e_{sq}+e_{ss}+e_{pm}-e_{pp}+e_{mm})=e_{pq}$. It implies that $\{e_{ps}+\frac{1}{2}e_{pp}+e_{mp}-\frac{1}{2}e_{mm}, e_{sq}+e_{ss}+e_{pm}-e_{pp}+e_{mm}\}=w$. As $\{\cdot, \cdot\}$ is symmetric, by (25), (27), (31) and (59), if $m\neq p\;\mathrm{and} \;m\neq q$ this yields
\begin{equation}
\{e_{mp}, e_{pm}\}=\frac {1}{2}\{e_{pp}, e_{pp}\}+\frac {1}{2}\{e_{mm}, e_{mm}\}=\{e_{pm}, e_{mp}\}.
\end{equation}
Choosing $s\neq q\; \mathrm {and} \;s\neq m$, then if $m\neq p\;\mathrm{and} \;m\neq q$ we have $(e_{ps}+e_{mq}+\frac {1}{2}e_{qq}-\frac{1}{2}e_{ss}-\frac {1}{2}e_{mm})\circ (e_{sq}+e_{qm}-e_{qq}+e_{mm})=e_{pq}$. It is clear that $\{e_{ps}+e_{mq}+\frac {1}{2}e_{qq}-\frac{1}{2}e_{ss}-\frac {1}{2}e_{mm}, e_{sq}+e_{qm}-e_{qq}+e_{mm}\}=w$. Because $\{\cdot, \cdot\}$ is symmetric, if $m\neq p\;\mathrm{and} \;m\neq q$ it follows from (25), (27), (31) and (57) that
\begin{equation}
\{e_{mq}, e_{qm}\}=\frac {1}{2}\{e_{mm}, e_{mm}\}+\frac{1}{2}\{e_{qq},e_{qq}\}=\{e_{qm},e_{mq}\}.
\end{equation}
Since $(e_{pp}+\frac{5}{4}e_{pq}+e_{qp}+e_{qq})\circ (e_{pp}-\frac{3}{4}e_{pq}-e_{qp}+e_{qq})=e_{pq}$, then we have that $\{e_{pp}+\frac{5}{4}e_{pq}+e_{qp}+e_{qq}, e_{pp}-\frac{3}{4}e_{pq}-e_{qp}+e_{qq}\}=w$. As $\{\cdot, \cdot\}$ is symmetric, using (25), (27) and (38), it can be reduced to
\begin{equation}
\{e_{pq}, e_{qp}\}=\frac{1}{2}\{e_{pp}, e_{pp}\}+\frac{1}{2}\{e_{qq}, e_{qq}\}=\{e_{qp}, e_{pq}\}.
\end{equation}

Step 10. If $j\neq p\; \mathrm {and} \;j\neq q$, we can verify that $(e_{pp}+e_{pq}-4e_{qp}-3e_{qq}+e_{qj}+e_{jj})\circ (2e_{pp}-4e_{qp}-e_{qq}-4e_{jp}-2e_{jq}+e_{jj})=e_{pq}$. It follows that $\{e_{pp}+e_{pq}-4e_{qp}-3e_{qq}+e_{qj}+e_{jj}, 2e_{pp}-4e_{qp}-e_{qq}-4e_{jp}-2e_{jq}+e_{jj}\}=w$. Since $\{\cdot, \cdot\}$ is symmetric, applying (25), (27), (31), (38), (39), (40), (53), (54), (57), (59), (61) and (62), if $j\neq p\; \mathrm {and} \;j\neq q$ we have that
\begin{equation}
\{e_{qj}, e_{jp}\}=\{e_{qp}, e_{qq}\}=\{e_{qp}, e_{pp}\}=\{e_{qq}, e_{qp}\}=\{e_{pp}, e_{qp}\}=\{e_{jp}, e_{qj}\}.
\end{equation}

Step 11. Now, if $i\neq m\; \mathrm {and} \;j\neq k$ we can conclude from step 1, step 4, step 5 and step 6 that
\begin{equation}
\{e_{ij}, e_{km}\}=0.
\end{equation}
If $i,j,m$ are distinct, it follows from step 2, step 7, step 8 and step 10 that
\begin{equation}
\{e_{ij}, e_{jm}\}=\{e_{im}, e_{mm}\}=\{e_{ii}, e_{im}\}=\{e_{im}, e_{ii}\}=\{e_{mm}, e_{im}\}=\{e_{jm}, e_{ij}\}.
\end{equation}
If $i\neq m$ and $n=3$, from (60), (61) and (62) we have that
\begin{equation}
\{e_{im}, e_{mi}\}=\frac{1}{2}\{e_{ii}, e_{ii}\}+\frac{1}{2}\{e_{mm}, e_{mm}\}.
\end{equation}
If $i\neq m$ and $n>3$, from step 3 and step 9 we have the same equations as well. So we complete the proof. $\Box$

\section*{3. Several applications}
~

Now, we will give two applications of the theorems above.\\~\\{\bf Definition 3.1.} We say that $G\in M_n(R)$ is a Jordan all-multiplicative point in $M_n(R)$ if for every $M_n(R)$-module $X$ and every Jordan multiplicative $R$-linear map $\varphi$ : $M_n(R)\to X$ at $G$ (i.e. $\varphi (S\circ T)=\varphi (S)\circ \varphi (T)$ for any $S,T\in M_n(R),\: S\circ T=G$) with $\varphi (I)=I$ is a multiplicative mapping in $M_n(R)$.\\~\\{\bf Corollary 3.2.} Every matrix units $e_{ij}$ in $M_n(R)$, $n\geq 3$, is a Jordan all-multiplicative point.\\~\\{\bf Proof.} Let $X$ be an $M_n(R)$-module, $\varphi$ be a Jordan multiplicative $R$-linear map at $e_{ij}$. Then it follows
\begin{equation}\nonumber
\varphi (I)=I, \;\varphi (S\circ T)=\varphi (S)\circ \varphi (T)\ \mathrm{for\;all}\: S,T\in M_n(R),\:S\circ T=e_{ij}.
\end{equation}
Set $\{S,T\}=\varphi (S) \circ \varphi (T)$ for any $S,T\in M_n(R)$, thus $\{\cdot,\cdot\}$ is a symmetric $R$-bilinear map. Since $\varphi$ is a multiplicative map at $e_{ij}$, we have
\begin{equation}\nonumber
\{S,T\}=\varphi (S)\circ \varphi (T)=\varphi (e_{ij}) \ \mathrm{for\;all}\: S,T\in M_n(R)\;\mathrm{with}\;S\circ T=e_{ij}.
\end{equation}
As $e_{ij}$ is a Jordan product determined point in $M_n(R)$, then there exists an $R$-linear map $\phi$ : $M_n(R)^2\to X$ such that $\{S,T\}=\phi(S\circ T)$ for all $S,T\in M_n(R)$. So
\begin{equation}\nonumber
\varphi (S)\circ \varphi (T)=\{S,T\}=\phi (S\circ T),\; \forall \:S,T\in M_n(R).
\end{equation}
Set $S=I$ in the equation above, then we have $\varphi (T)=\phi (T)$ for every $T\in M_n(R)$. It follows that
\begin{equation}\nonumber
\varphi(S\circ T)=\varphi (S)\circ \varphi (T),\;\forall\; S,T\in M_n(R).
\end{equation}
Hence $\varphi$ is a multiplicative mapping in $M_n(R)$. The proof is completed. $\Box$ \\~\\{\bf Definition 3.3.} We say that $H\in M_n(R)$ is a Jordan all-derivable point in $M_n(R)$ if for every $M_n(R)$-module $X$ and every Jordan  derivable $R$-linear map $\varphi$ : $M_n(R)\to X$ at $H$ (i.e. $\varphi (S\circ T)=\varphi (S)\circ T+S\circ \varphi (T)$ for any $S,T\in M_n(R),\: S\circ T=H$) with $\varphi (I)=0$ is a Jordan derivation in $M_n(R)$.\\~\\{\bf Corollary 3.4.} Every matrix units $e_{ij}$ in $M_n(R)$, $n\geq 3$, is a Jordan all-derivable point.\\~\\{\bf Proof.} Let $X$ be an $M_n(R)$-module, $\tau$ be a Jordan derivable $R$-linear map at $e_{ij}$. Then we have
\begin{equation}\nonumber
\tau (I)=0,\; \tau  (S\circ T)=\tau  (S)\circ T+S\circ \tau  (T)\;\mathrm{ for\; all}\;S,T\in M_n(R),\;S\circ T=e_{ij}.
\end{equation}
Set $\{S,T\}=\tau  (S)\circ T+S\circ \tau  (T)$ for any $S,T\in M_n(R)$, so $\{\cdot,\cdot\}$ is a symmetric $R$-bilinear map. Since $\tau$ is a derivable map at $e_{ij}$, it follows
\begin{equation}\nonumber
\{S,T\}=\tau  (S)\circ T+S\circ \tau  (T)=\tau (e_{ij})\;\mathrm{ for\; all}\;S,T\in M_n(R)\;\mathrm{with}\;S\circ T=e_{ij}.
\end{equation}
As $e_{ij}$ is a Jordan product determined point in $M_n(R)$, then there exists an $R$-linear map $\psi$ : $M_n(R)^2\to X$ such that $\{S,T\}=\psi(S\circ T)$ for all $S,T\in M_n(R)$. Hence
\begin{equation}\nonumber
\tau (S)\circ T+S \circ \tau (T)=\{S,T\}=\psi (S\circ T),\; \forall \:S,T\in M_n(R).
\end{equation}
Set $S=I$ in the equation above, then we have $\tau (T)=\psi (T)$ for every $T\in M_n(R)$. It follows that
\begin{equation}\nonumber
\tau(S\circ T)=\tau (S)\circ T+S\circ \tau (T),\; \forall \:S,T\in M_n(R).
\end{equation}
Then $\tau$ is a Jordan derivation in $M_n(R)$. The proof is completed. $\Box$\\

\end{document}